\documentclass{amsart}

\newcommand{\linf}{\ensuremath{\lim_{n \rightarrow \infty}}}
\newtheorem{thm}{Theorem}

\begin{document}

\begin{abstract}
The property that a 1-1 function from the set of natural numbers
, $\mathcal{N}$,
to itself
preserves the density of subsets of $\mathcal{N}$ is shown to be
equivalent to a condition on the covering of intervals in the range of the
function by images of intervals in the domain of the function.
\end{abstract}

\title{Density Preserving Functions}
\author{Paul J. Huizinga}
\address{ICL-1 L208\\
City College of San Francisco\\
50 Phelan Ave\\
San Francisco, CA 94112}
\email{phuizing@ccsf.cc.ca.us}
\maketitle

\section{Density}

If $\mathcal{N}$ is the set of natural numbers and $\mathcal{F}$ is a
finite subset of $\mathcal{N}$ , then $\| \mathcal{F} \|$ denotes the
number of elements of $\mathcal{F}$. For $\mathcal{S}$, an arbitrary
subset of $\mathcal{N}$, let $\mathcal{S}_n$ denote the set of elements
of $\mathcal{S}$ less than or equal to $n$. If the limit
\[d(\mathcal{S}) = \linf \frac{\| \mathcal{S}_n \|}{n} \]
exists, $\mathcal{S}$ is said to have the density $d(\mathcal{S})$.
The concept of density and variations on it 
occur in several areas of mathematics,  
for example, probability theory \cite[ch. VIII sec. 4]{feller}, 
algebraic number theory \cite[ch. VIII sec. 4]{lang}, 
and the number theoretic study of subsets of the natural numbers
 \cite[ch.V]{halb}.

An example of a set which does not have a density is

\[ \mathcal{S} = \{n |\; 2^{2m} \leq n < 2^{2m+1} 
,\; m = 0, 1, 2, \ldots \}.\]
There are $4^m$ elements of $\mathcal{S}$
associated with the value, m.
If $ n \; = \; 2^{2m + 1} -1$ and $ n' \; = \; 2^{2m + 2} - 1 $,
then $\| \mathcal{S}_n \| \; = \; \| \mathcal{S}_{n'} \| $ is
\[ \sum_{i=0}^m 4^i \; = \; 
  \frac{4^{m+1} - 1}{3}. \]
So that $ \| \mathcal{S}_n \| / n$ is
\[ \frac{ \frac{ 4^{m+1} -1 }{ 3 } }{ 2^{2m + 1} - 1 }
 \; = \; \frac{1}{3} \left( \frac{ 2^{2m +2} -1 }{ 2^{2m+1} -1} \right)
 \; = \; \frac{2}{3} \left( \frac{ 2^{2m+1} -\frac{1}{2} }{ 2^{2m+1} -1}
\right) 
\]
and $ \| \mathcal{S}_{n'} \| / n'$ is
\[ \frac{ \frac{ 4^{m+1} -1 }{3} }{ 2^{2m+2} -1 } 
 \; = \; \frac{1}{3} \left( \frac{ 2^{2m+2} - 1 }{ 2^{2m+2} - 1} \right).
\] 
Therefore $\| \mathcal{S}_n \| / n$ has a $\limsup$ of 2/3 and a 
$\liminf$ of 1/3. 

An interval in the set of natural numbers 
is a sub-set of $\mathcal{N}$ of the form: 
\[ I = [a, b] = \{ n | \; a \leq n \leq b \}. \]
For such an interval I, $\mu (I)$ is defined to be b/a.
The interval, $I=[a,b]$, is said to be an m-interval
( $ m > 1,\; M \in \Re $ -- the set of real numbers), 
if $\frac{b}{a} \leq m < \frac{b+1}{a}$ 
or equivalently $ma - 1  <  b  \leq ma.$
$I$ is said to be a $\;$ +m-interval (plus-m-interval) 
, if $ \mu(I) > m$. 

If $I = [a,b]$, call 
\[ \frac{ \| \mathcal{S} \cap I \| }{ \| I \| } 
  \; = \; \frac{ \| \mathcal{S} \cap I \| }{b-a+1} \]
the density of $\mathcal{S}$ in $I$.
The intervals, $ [ 2^{2m}, 2^{2m +1} - 1 ] $
contained in the set with no density, $\mathcal{S}$,
in the example above
and the intervals, $ [ 2^{2m+1}, 2^{2m+2} - 1 ]$
contained in its complement have a
$\mu > 1.5$ for $m \geq 1$.
Therefore, if $ I_x $ is the 1.5-interval
which has x as left endpoint,
the density of $\mathcal{S}$ in $I_x$ does not converge
(in fact oscillates between 1 and 0) 
as $x$ goes to infinity.
That this is a characteristic of sets which fail to have a density is
shown by the following:

\begin{thm}
For any set, $ \mathcal{S}, \; d(\mathcal{S})$ exists and is equal to $D$
iff for any $\epsilon > 0$ 
($\epsilon \in \Re$)
and $m > 1 \; ( \; m \in \Re )$ there is an $N \geq 1 
  \; ( N \in \mathcal{N} )$
 such that for any +m-interval $I = [a,b]$ with $a > N$
\[ \left| \frac{\| \mathcal{S} \cap I \|}{ \| I \| } - D \right| 
  < \epsilon . \]
\end{thm}

First, suppose $d(\mathcal{S}) = D$. Given $\epsilon$  and $m$, 
let $0 < \epsilon ' < \frac{m-1}{m+1}\epsilon$.
Since $d(\mathcal{S}) = D$ there is an $N$ such that, if $n>N$
\[ \left| \frac{\| \mathcal{S}_n \| }{n}  - D \right| < \epsilon ' .\]
If $I = [a,b]$ with $a>N+1$ and $\mu(I) > m$, then 
\[ Db - \epsilon'b < \|\mathcal{S}_b \| < Db + \epsilon'b \]
\[ -D(a-1) - \epsilon'(a-1) < - \|\mathcal{S}_{a-1}\|
 < -D(a-1) + \epsilon'(a-1) . \]
Since $ \| \mathcal{S} \cap I \|  =  \| \mathcal{S}_b \| -
\| \mathcal{S}_{a-1} \|$, we have
\[ D(b-(a-1) - (b+(a-1))\epsilon' < \| \mathcal{S} \cap I \| <
D(b-(a-1)) + (b+(a-1))\epsilon' . \]
When $c>0$, $\frac{x-c}{x+c}$ is an increasing function of $x$ for 
$x \neq -c$. So that $b > ma$ implies 
\[ \frac{b-(a-1)}{b+(a-1)} \; > \; \frac{ma - a + 1}{ma + a - 1}
 \; > \; \frac{ma-a}{ma+a} \; = \; \frac{m-1}{m+1} \]
and
\[ (b-(a-1))\epsilon \; = \; (b+(a-1)) \frac{b-(a-1)}{b+(a-1)} \epsilon
\]
\[ \; > \; (b+(a-1)) \frac{m-1}{m+1} \epsilon \; > \; (b+(a-1)) \epsilon'
\]
(note the prime on the final $\epsilon$) giving
\[ D(b-(a-1)) - \epsilon (b-(a-1)) \; < \; D(b-(a-1))
   - \epsilon' (b+(a-1)) \] 
\[ \; < \: \| \mathcal{S} \cap I \| \]
\[ \; < \; D(b-(a-1)) + \epsilon'(b+(a-1))
   \; < \; D(b-(a-1)) + \epsilon (b-(a-1))
\]
or
\[ \left| \frac{ \| \mathcal{S} \cap I \|}{ \| I \| } 
   - D \right| \; = \;
   \left| \frac{ \| \mathcal{S} \cap I \|}{b-(a-1)} - D \right|
   \; < \;   \epsilon. \] 

Now suppose the second condition is met, that is the density of
$\mathcal{S}$ in +m-intervals approaches $D$ asymptotically for
for any $m>1$.
Given $\epsilon$, choose $\epsilon' < \epsilon / 3$
and $m > 3 / \epsilon$. 
If $N$ is the value given by the second condition and $I = [a,b]$ is
an interval with $a>N$ and $\mu(I) > m = 3 / \epsilon$, then
\[ D(b-(a-1)) - \epsilon' (b-(a-1)) \; < \; \| \mathcal{S} \cap I \|
\; < \; D(b-(a-1)) + \epsilon' (b-(a-1)) \]
also
\[ \| \mathcal{S} \cap I\| \; < \; \|\mathcal{S}_b \|
\; < \; \| \mathcal{S} \cap I \| + a. \]
Now
\[ \| \mathcal{S} \cap I \| \: > \; D(b-(a-1)) - \epsilon' (b-(a-1))
\; = \; Db - \epsilon'b - D(a-1) + \epsilon'(a-1) \] 
\[ > \; Db - \epsilon'b - D(a-1) \]
\[ > \; Db - \frac{\epsilon}{3}b - \frac{\epsilon}{3}b \]
because $\epsilon' < \epsilon / 3$, $a < (\epsilon / 3)b$,
and $D \leq 1$. So that
\[ \| \mathcal{S}_b \| \; > \;  \| \mathcal{S} \cap I \|
 \; > \; Db - \epsilon b. \]

On the other hand
\[ \| \mathcal{S}_b \| \; < \; \| \mathcal{S} \cap I \| + a
\; < \; Db + \epsilon' b + a - D(a-1) - \epsilon' (a-1) \]
\[ < \; Db + \epsilon' b +a \; < \; Db + \frac{\epsilon}{3} b
+ \frac{\epsilon}{3} b. \]
So that
\[ \| \mathcal{S}_b \| \; < \; Db + \epsilon b. \]
Therefore, for any $b > (3 / \epsilon) N$, 
\[ \left| \frac{ \| \mathcal{S}_b \| }{b} -D \right| \; < \; \epsilon \]
and $d(\mathcal{S}) = D$.

\section{Density Preserving Functions}

A function, $f: \mathcal{N} \rightarrow \mathcal{N},$ 
is said to preserve density if it is one to one 
and whenever $d( \mathcal{S} ) = D$, we have 
$ d(f( \mathcal{S})) = D$.
The next theorem describes density preserving functions in terms of
the existence for any $p>1$ of a value, $m>1$, such that (asymptotically)
any collection of m-intervals whose images cover a +p-interval
will have a sub-collection, $\mathcal{C}$,  
with a specified goodness-of fit.   
The goodness-of fit is given by two values, $q$ and $r$.
The value, $q$ (the inclusion factor), in the theorem, 
is the fraction of the covering set which is in the
covered interval and may be thought of as close to 1. The value,
$r$ (the omission factor), is the fraction of the covered interval not in 
the union of the images of the sub-collection, $\mathcal{C}$,
and may be thought of as close to 0.

\newpage
\begin{thm}
Let $f:\mathcal{N}\rightarrow\mathcal{N}$ be 1-1, then $f$ preserves
density
iff
\[ \forall p \in \Re, \; p>1 \]
\[ \forall q \in \Re, \; 0<q<1 \]
\[ \forall r \in \Re, \; 0<r<1 \]
\[ \exists m \in \Re, \; m>1 \]
\[ \exists N \in \mathcal{N}, \; N \geq 1 \]
such that if $I = [a,b]$ is a +p-interval with $a>N$,  
$\{J_i | \, i=1 \ldots k \}$ is any disjoint collection of m-intervals
with
\[ I \subset \cup_{i=1}^k f( J_i ), \]
and $ \mathcal{C} \; = \;
 \{ J_i | \;  \| I \cap f (J_i) \| \:\geq \; q \| J_i \| \}; $
then for 
$ T \; = \; \cup \{ J_i | \; J_i \notin \mathcal{C} \} $
we have
\[ \| f(T) \cap I \| \; < \; r \| I \|. \]
\end{thm}

First suppose the covering condition holds and $d(\mathcal{S})=D.$
If a $p>1$ and an $\epsilon > 0$ are given, 
we may assume without loss of generality
that $\epsilon < 3D$ to simplify the choice of $r$ below.
It will be shown that there is an $N''$ 
such that if $I = [a,b]$, $\mu(I) > p$, and $a>N''$; then
\[ \left| \frac{ \| f(\mathcal{S} ) \cap I \| }{ \| I \| }
 \; -\; D \right| \; < \: \epsilon . \]

To do this, choose
\[ 0 \; < \; r \; < \; min \left( \frac{\epsilon}{3},  
\frac{ \frac{\epsilon}{6}}{ D - \frac{\epsilon}{3}} \right) \]
\[  max \left( \frac{1}{1 + \frac{\epsilon}{3}},
\frac{D + \frac{\epsilon}{3}}{ D + \frac{\epsilon}{2}} \right)
\; < \; q \; < \; 1 . \]

Let $m$ and $N$ be the values given by the hypothesis for $p$ , $q$, and
$r$. Since $d(\mathcal{S}) = D$,
by theorem 1, there is an $N'$ such that
for $J = [c,d]$, $c>N'$, $J$ an m-interval, then
\[ \left| \frac{ \| J \cap \mathcal{S} \| }{ \| J \| } - D \right|
\; < \: \frac{\epsilon}{3}. \]
Choose $N'' \; > \: max \{ f(i) | \; i \leq m N' \}$ and also 
$N'' \; > \; N$.

Let $I=[a,b]$ be a +p-interval with $a > N''$
and $\{ J_1 \ldots J_k \}$ be a
disjoint collection of m-intervals such that
$ I \; \subset \: \cup_{i=1}^k f( J_i ). $

If $J_i = [c_i,d_i]$ with $c_i \leq N'$, 
then $d_i \, \leq \, mc_i \, \leq \, mN'$.
So that $f( J_i ) \cap I = \emptyset$
and $J_i$ is not in $\mathcal{C} \; = \;
 \{ J_i | \;  \| I \cap f (J_i) \| \:\geq \; q \| J_i \| \}$.
If $J_i = [c_i, d_i]$ is in $\mathcal{C}$,
then $c_i > N'$ and
\[ ( D - \frac{ \epsilon }{3} ) \| J_i \|
 \; < \; \| \mathcal{S} \cap J_i \|
 \; < \; ( D + \frac{ \epsilon }{3} ) \| J_i \| . \]
Let
\[ K \, = \, \cup \{ f(J_i) | \; J_i \in\mathcal{C} \}. \]
Set $I_1 = I \cap K, \; I_2 = I - I_1$.
By hypothesis, $\| I_2 \| < r \| I \| < \frac{\epsilon}{3} \| I \|$.
By the definition of $K$, $ q \| K \| \leq \| I_1 \|$,
$ \| K \| \leq \frac{1}{q} \| I_1 \|$.
So that
\[ \|K - I_1 \| \; = \; \| K \| - \| I_1 \| \; \leq \; \frac{1-q}{q}
\| I_1 \|. \]
Since
\[ q \; > \; \frac{1}{1+ \frac{ \epsilon }{3} }, \]
we have
\[ \frac{1-q}{q} \: < \: \frac{ \epsilon }{3} \]
and since $ \| I_1 \| \leq \| I \| $,
\[ \| K - I_1 \| \; < \; \frac{ \epsilon }{3} \| I_1 \| 
\; \leq \; \frac{ \epsilon }{3} \| I \|. \]

The function, $f$, is 1-1,
so for $\mathcal{S}$ the set of density $D$
\[ \| f( \mathcal{S} ) \cap K \| \; > \; 
\left( D - \frac{ \epsilon }{3} \right) \| K \| \; > \:
\left( D - \frac{ \epsilon}{3} \right) \| I_1 \| \; > \;
\left( D - \frac{ \epsilon}{3} \right) (1-r) \| I \| . \]
By the choice of $r$
\[ r \; < \; \frac{ \frac{ \epsilon }{6} }{ D - \frac{ \epsilon }{3} } \]
and
\[ 1-r \; > \; \frac{ D - \frac{ \epsilon }{2} }{
D - \frac{ \epsilon }{3} }. \]
So that
\[ \| f( \mathcal{S} ) \cap K \| \; > \left( D - \frac{ \epsilon }{2}
\right) \| I \| \]
and
\[ \| f( \mathcal{S} ) \cap I \| \; \geq \: \| f( \mathcal{S} ) \cap I_1
\|
\; \geq \; \| f( \mathcal{S} ) \cap K \| \, - \, \| K - I_1 \| \]
\[ > \; \left( D - \frac{ \epsilon }{2} \right) \| I \| \, - \,
\frac{ \epsilon }{3} \| I \| \] 
\[ \; > \; ( D - \epsilon ) \| I \| . \]

On the other hand, we have
\[ q \; > \; \frac{ D + \frac{ \epsilon }{3} }{ 
D + \frac{ \epsilon }{2} } \]
or
\[ \frac{1}{q} \; < \; \frac{ D + \frac{ \epsilon }{2} }{
D + \frac{ \epsilon }{3} } . \]
So that
\[  \| f( \mathcal{S} ) \cap I \|
 \; \leq \; \| f( \mathcal{S} ) \cap K \| + \| I_2 \|
 \; < \;  \left( D + \frac{ \epsilon }{3} \right) \| K \| + \| I_2 \| \]
\[ < \; \left( D + \frac{ \epsilon }{3} \right) \frac{1}{q} \| I \|
 + \frac{ \epsilon}{3} \| I \|
 \; < \; \left( D + \frac{ \epsilon }{2} \right) \| I \|
 + \frac{ \epsilon}{3}  \| I \| \]
\[ < \; ( D + \epsilon ) \| I \|. \]

Combining the two inequalities 
\[ (D - \epsilon )\| I \| \; < \; \| f( \mathcal{S} ) \cap I \|
\; < \; (D + \epsilon ) \| I \|  \]
or
\[ \left| \frac{ \| f( \mathcal{S} ) \cap I \| }{\| I \|} - D \right|
 \; < \; \epsilon. \]
Since this holds for any value of $p>1$, by theorem 1 we have that
$d( f( \mathcal{S} ) ) = D.$

In the other direction, suppose that the covering condition fails to
hold.  Then there exist $p$, $q$, and $r$; such that for all $m$ and all
$N$ there is a +p-interval 
\linebreak
$I = [a,b]$ with $a > N$ and a collection of
disjoint m-intervals
 \[ \{ J_i | i = 1 \ldots k \} \]
such that $I$ is contained in the union of the images of the $J_i$
and for
\[ T \; = \; \cup \{ J_i | \; \| I \cap f (J_i) \| \:< \; q \| J_i \| \}
,\]
we have
\[ \| f(T) \cap I \| \; \geq \; r \| I \| . \]
Since this is also true for $ r' < r$, we may assume that $ r <
  \frac{1}{2}.$

The idea is to construct a set whose density is not preserved.
Let $\lfloor x \rfloor =$ the greatest integer $\leq x.$
Suppose a set, $\mathcal{S}$, of density, D, is being constructed.
If, at some stage in the construction, $ \lfloor Dn \rfloor $ values 
less than or equal to $n$
have been included in $\mathcal{S}$ and all other values
less than or equal to $n$ have been excluded.
Then
\[ D - \frac{1}{n} \; < \; \frac{\| S_n \| }{ n }  \; \leq \: D. \]
If there are no constraints on the choice of elements of $\mathcal{S}$,
$i$ can be
chosen to be in $\mathcal{S}$ whenever $ \lfloor D(i-1) \rfloor \, < \, 
\lfloor Di \rfloor $ and the above inequality will be true for every $n$.

To construct a sequence whose density is not preserved,
some constraints must be placed on the choice of elements of
$\mathcal{S}$.
At the k-th stage of the construction, these constraints will consist of
choosing certain elements of $( 1 + \frac{1}{k} )$-intervals.
If $J = [c,d]$ is such a $( 1 + \frac{1}{k} )$-interval,
it will be the case that
$c$ is greater than $4k$ and
$ \lfloor D(c-1) \rfloor $ values less than $c$
will have been assigned to $S$. 
When the construction reaches $d$,
$ \lfloor Dd \rfloor $
elements of $\mathcal{S}$ will have been chosen.
Therefore, no more than
$\lfloor D \| J \| \rfloor + 1$
elements of $J$ will have been added to $\mathcal{S}$.
The density of $\mathcal{S}$ in $J$ will be less than or equal to
$ D \; + \; \frac{1}{ \| J \| }. $ 

Since, $J$ is a $(1 + \frac{1}{k})$-interval,
\[ (1 + \frac{1}{k} )c - 1 \; < \; d \; \leq \; ( 1 + \frac{1}{k} )c. \]
So that,
\[ \frac{c}{k}  \; < \; d - (c-1) \, = \, \| J \|
 \; \leq \; \frac{c}{k} + 1. \] 
Therefore, even if all of the elements of $J$ are added to $\mathcal{S}$
we will have
\[ \frac{ \| \mathcal{S}_d \| }{ d} \; \leq \;
 \frac{ \| \mathcal{S}_{c-1} \| \, + \, \| J \| }{ d }
\; < \; 
\frac{ D(c-1) + \frac{c}{k} + 1 }{ c + \frac{c}{k} - 1 } \]
\[ = \; \frac{ D( c + \frac{c}{k} - 1 ) 
  \, + \, (1-D) \frac{c}{k} + 1}{ c + \frac{c}{k} - 1 } \]
\[ = \; D \, + \, \frac{ (1-D) \frac{c}{k} + 1}
  { c + \frac{c}{k} - 1 } \]
\[ = \; D \, + \, \frac{ (1-D) + \frac{k}{c} }
  { k + 1 - \frac{k}{c} } \]
\[ < \; D \, + \, \frac{2}{k}. \]
The last inequality holds because $ 0 \leq D \leq 1$ and $ c >4k.$
A similar argument holds if no elements of $J$ are added to $\mathcal{S}$.
Therefore, for any $n$, $c \leq n \leq d$
it will be true that: 
\[\left| \frac{ \| \mathcal{S}_n\| }{ n } \, - \, D \right| \; < \;
\frac{2}{k} \]
and $d(\mathcal{S})$ will exist and be equal to $D$.

Given a function, $f$, for which the covering condition fails to hold with
values $p$, $q$, and $r$;
a set, $\mathcal{S}$, of density $D = \frac{1-q}{2}$ will be constructed.
At the end of stage $k$ all values less than or equal to $L_k$ will have
been included in or excluded from the set $\mathcal{S}$
and the membership of values above $L_k$ will be undetermined.
Set $L_0 = 0$.

At stage $k$, set 
\[ M_k = \max \left( ( 1 + \frac{1}{k} ) L_{k-1}, \; 
 \frac{4k (1-r)}{r(1-q)} \right) 
\]
and $N_k = \max\{ f(x) | \; x \leq M_k \} + 1$
Then for $(1+\frac{1}{k})$-intervals,
$[c,d]$ with $c > M_k$ 
(since $r$ is -- by assumption -- less than $ \frac{1}{2}$)
\[ d-(c-1) \; > \; \frac{1}{k}c  
  \; > \; \frac{1}{k} M_k  
  \; > \; \frac{ 4 (1-r) }{ r (1-q) }
  \; > \; \frac{4}{1-q} \]
and $ (1-q)(d-(c-1)) \; > \; 4 $.
This means that
\[ 2 \; < \; \frac{1-q}{2} (d-(c-1)) \] 
\[ \frac{1-q}{2} (d-(c-1)) + 2 \; < \; 
   (1-q)(d-(c-1)) \]
so that there will be no problem with choosing
 $ \lfloor \frac{1-q}{2}(d-(c-1)) \rfloor + 1$
elements out of a subset of $[c,d]$ containing
at least $(1-q)(d-(c-1))$ elements.
Also,
\[ c \; > \; \frac{4k}{1-q} \; > \; 4k \]
as mentioned above.

The fact that the covering condition does not hold implies that
for $N = N_k$ and $m = (1 + \frac{1}{k})$ there is a +p-interval,
$I = [a,b]$ with $a>N_k$ and a collection, 
$\{ J_i\}$ of disjoint $(1 + \frac{1}{k})$-intervals 
whose images cover $I$ such that the
$f(J)$'s with inclusion factor less than $q$ 
contain more than $r \| I \|$ elements of $I$.
Since $N_k > \max \{ f(x) | x \leq M_k \}$ and 
$M_k \geq ( 1 + \frac{1}{k} ) L_{k-1}$,
no value in a $( 1 + \frac{1}{k} )$-interval whose image intersects
I has been included in or excluded from $\mathcal{S}$ 
at the end of stage $ k - 1$. 

Starting at $L_{k-1} + 1$ the k-th stage of the construction proceeds in
ascending order.
If all values less than $x$ have been assigned to $\mathcal{S}$ 
or $ \neg \mathcal{S}$ and
$x$ is not in a $J$ whose image intersects $I$, 
then $x$ is assigned 
to $\mathcal{S}$ if and only if
$ \lfloor D(x-1) \rfloor \, < \, \lfloor Dx \rfloor $.
When an interval, $J = [c,d]$, in the given collection 
whose image intersects $I$ is reached,
calculate how many elements of $J$ must be added to $\mathcal{S}$ in order
for $ \| \mathcal{S}_{d} \| \, = \, \lfloor Dd \rfloor $.
At most $ \lfloor D(d-(c-1)) \rfloor + 1 $ will be needed.
Choose as many as possible of them from the elements of $J$ whose
images are not in $I$.
In the case of the intervals not in $\mathcal{C}$, all of the
elements can be chosen so that their image is not in $I$.
In the other intervals,
since an element whose image is in $I$ is included in $\mathcal{S}$
only if all elements whose images are not in $I$ have been included, 
the proportion of elements in $J \cap f^{-1}(I)$
that are assigned to $\mathcal{S}$ is less than or equal to 
 $D \; + \; \frac{1}{ \| J \| }$.
When the construction has assigned all the elements of the $J$-s,
set $L_k$ equal to the last value considered.

This choice of elements of $\mathcal{S}$ yields
\[ \| f(\mathcal{S}) \cap I \| \; \leq \; \left[
  \sum_{J \in \mathcal{C} } \left( D + \frac{1}{ \| J \| } \right) 
  \| f(J) \cap I \|  \right]
  \; + \; \left[ 0 \cdot \sum_{J' \notin \mathcal{C} } \| f(J') \| \right]
 \]

Those $J = [c,d]$ whose images intersect $I$ have 
\[ c \; > \; \frac{ 4k (1-r) }{ r (1-q) } \]
Which means
\[ \| J \| \; > \; \frac{c}{k} \; > \; \frac{ 4 (1-r) }{ r (1-q)} \]
or
\[ \frac{1}{ \| J \| } \; < \; \frac{ r (1-q) }{ 4 (1-r) }
  \; = \; \left( \frac{1}{ 1-r } \right) \frac{r}{2} 
   \left( \frac{ 1-q }{ 2 } \right)
   \; = \; \frac{ \frac{r}{2} D }{ 1-r } \]
and
\[ D \; + \frac{1}{ \| J \| } 
  \; < \: \frac{ ( 1 - \frac{r}{2} ) D }{ 1-r } . \]
 
So that
\[ \| f( \mathcal{S} ) \cap I \|
  \; < \;\left( \frac{ ( 1 - \frac{r}{2} ) D }{ 1-r } \right) 
  \cdot \sum_{J \in \mathcal{C} } \| f(J) \cap I \|
  \; \leq \; \left( \frac{ ( 1 - \frac{r}{2} ) D }{ 1-r } \right)
  \cdot (1-r) \| I \| \]
and 
\[ \frac{ \| f(\mathcal{S}) \cap I \| }{ \| I \| } 
  \; < \;D - D \frac{r}{2} . \]

When the construction is completed, for any $N$, we have a 
+p-interval whose elements are greater than $N$ 
and whose local density is  at least $D\frac{r}{2}$ 
less than than $D$.
Therefore the density of $f(\mathcal{S})$ is not $D$ and $f$ does not
preserve density.

\section{ An Example: The $2^n$ Shuffle}

The $2^n$ shuffle, $sh()$, is defined as follows:
\[ sh(k) = \begin{cases} k &: \quad k < 4 \\
   2^i + 2j &: \quad k = 2^i +j, \quad i > 2, \quad 0 \leq j < 2^{i-1} \\
   2^i + 2j +1 &: \quad k = 2^i + 2^{i-1} + j, \quad i > 2, \quad 
   0 \leq j < 2^{i-1}.\\
\end{cases} \]
Informally, $sh()$ shuffles the numbers in $[ 2^i, \; 2^{i+1} -1 ]$ 
for $i \geq 2$
and its inverse, $sh^{-1}()$, deals the even numbers in that interval
to the lower half of the interval and the odd numbers to the upper half. 

Since $sh()$ is 1-1 and onto, when applying theorem 2, we can work in
the domain of $sh()$ as easily as in the range.
That is to say, we can consider coverings of the inverse image of a
+p-interval in the range by m-intervals in the domain.

If a +p-interval, $I$, contains all of $[ 2^i, 2^{i+1} - 1 ]$, 
its inverse image, $sh^{-1}(I)$, will also contain that interval.
If $I$ contains more than one but less than $ 2^i -1 $  
of the members of $[ 2^i, 2^{i+1} - 1 ]$, the inverse image of the
intersection of $I$ with that interval will consist of two intervals --
the even numbers going to the lower interval and the odd to the upper.
Therefore the inverse image of an interval under $sh()$ will consist of at
most 3 intervals, the even numbers being dealt to a lower interval at one
end and the odd to a higher at the other.

Next, consider the covering of a +p-interval, $I = [a, b]$,
 by a disjoint collection, $\{ J_i \}$, of m-intervals.
Assume that $ m \leq \sqrt[3]{p}$, so that at least one of the $J_i$
is completely contained in $I$.
The only $J_i$-s not entirely in $I$ or entirely in the complement of $I$
are the ones containing $a$ and $b$.
An m-interval containing a has at most $(m - 1)a + 1$ elements and if it 
intersects the complement of $I$, at most $(m - 1)a$ of them will be in
$I$.
A similar argument shows that there will be at most
$(m - 1)b$ elements in the intersection of $I$ and an
m-interval
containing $b$, but not entirely contained in $I$.
Let $\mathcal{C}'$ be the collection of $J_i$-s entirely contained in $I$.
$\mathcal{C}'$ is a sub-collection of 
$ \{ J_i | \; \| J_i \cap sh^{-1}(i) \| > q \| J_i \| \}$
for any $q < 1$, therefore if the omission factor for $\mathcal{C}'$
is $< r$, this will also be true for any inclusion factor, $q < 1$.
There are at most 
\[   (m - 1) ( b + a ) \]
elements in $I - \cup \mathcal{C}'$.
Since $I$ has $b-a+1$ elements,
the fraction of elements of $I$ not in $ \cup \mathcal{C}'$
is less than
\[ \frac{ (m-1) (b+a) }{ b-a }. \]
As in the proof of theorem 1, 
\[ \frac{b-a}{b+a} \; > \; \frac{pa - a }{ pa + a }
   \; = \; \frac{p-1}{p+1}. \] 
So that if m is close enough to 1,
\[0 \; < \; m-1 \; < \; \frac{p-1}{p+1} r
  \;\; \left( \; < \; \frac{b-a}{b+a} r \right), \]
then
\[ (m-1) \left( \frac{b+a}{b-a} \right) \; < \; r \]
and the omission factor of the sub-collection $ \mathcal{C}' $
is less than r.

Now let $I = [a, b]$ be a +p-interval whose inverse image, $sh^{-1}(I)$, 
is to be covered with an omission factor of $r$.
Since theorem 2 involves only the asymptotic properties of intervals,
we may require that $a$ be greater than a given  value
to be determined later.
We have
\[ \frac{b}{a} \; > \; p, \quad \frac{b}{p} \; > \; a,
   \quad -a \; > \; \frac{-b}{p}. 
\]
So that
\[ b-a+1 \; > \; b-a \; > b - \frac{b}{p} \; = \; \frac{p-1}{p}b. \]
That is $I$ has more than $\frac{p-1}{p}b$ elements and a sub-collection,
$ \mathcal{C}' $,
of a covering of disjoint intervals will have an omission factor less
than
$r$ if 
\[ \| sh^{-1}(I) - \cup \mathcal{C}' \| \; < \;
   \frac{p-1}{p} br. \]

$sh^{-1}(I)$ consists of at most 3 intervals.  The strategy will be to
discard intervals of sufficiently small $\mu$ and use covering intervals
whose $\mu$ is small enough that the sub-collection of intervals
contained in the inverse image will have an omission factor of 
less than $\frac{r}{3}$.

Exercising the option mentioned earlier, require $I = [a, b]$ to have
\[ a \; > \; \frac{6}{(p-1)r}. \]
Then, since $ a < \frac{b}{p}$,
\[ \frac{1}{6} \left( \frac{p-1}{p} br \right)
   \; > \; \frac{1}{6} (p-1)ar \; > \; 1 \]
and
\[ \frac{1}{3} \left( \frac{p-1}{p} br \right)
   \; > \; \frac{1}{6} \left( \frac{p-1}{p} br \right) + 1. \]

An interval with $k+1$ elements, $[x, x+k]$ has a $\mu$ of
\[ \frac{x + k}{x} \; = \; 1 \; + \; \frac{k}{x} \]
which is a decreasing function of $x$.
Therefore, the smallest $\mu$ for a component interval of $sh^{-1}(I)$
with $k+1$
elements will occur when this many odd elements are dealt upward
from the right hand side of $I$.

Let
\[ p' \; = \; 1 + \frac{1}{9} \left( \frac{p-1}{p} \right)r. \]
If $ 2^i + 1 < b <  2^{i+1} -2 $ and $ a < 2^i$,
then the right-most component of the inverse image of $ I = [a,b]$
will have the form (since $ 2^i + 2^{i-1} = \frac{3}{2} 2^i $ )
\[ \left[ \frac{3}{2} 2^i, \; \frac{3}{2} 2^i + k \right]. \]
If $\mu$ of this interval is less than or equal to $p'$, we have
\[ 1 + \frac{k}{ \frac{3}{2} 2^i } 
   \; \leq \; 1 + \frac{1}{9} \left( \frac{p-1}{p} \right)r
\]
\[ k \; \leq \; \frac{1}{9} \left( \frac{p-1}{p} \right)r ( \frac{3}{2} 2^i )
   \; < \; \frac{1}{6} \left( \frac{p-1}{p} \right) rb
\]
and
\[ k + 1 \; < \; \frac{1}{6} \left( \frac{p-1}{p} \right) rb + 1
  \; < \; \frac{1}{3} \left( \frac{p-1}{p} \right) br.
\]
Therefore any component interval of $sh^{-1}(I)$ with $ \mu \leq p'$ will have less than
\[ \frac{1}{3} \left( \frac{p-1}{p} \right) br \]
elements.

As shown earlier, a $+p'$-interval can be covered with an omission factor
of $ \frac{r}{3} $ by m-intervals where
\[ m \; < \; \sqrt[3]{p'} \]
and
\[ 0 \; < \; m-1 \; < \; \left( \frac{p'-1}{p'+1} \right) \frac{r}{3}
   \; = \; \frac{1}{3} \left( \frac{(p-1)r^2}{ 18p + (p-1)r } \right). 
\] 

 Up to two intervals in $sh^{-1}(I)$ of $\mu \leq p'$
can be ignored and $+p'$-intervals in $sh^{-1}(I)$ can be covered
with an omission factor of $\frac{r}{3}$ by m-intervals yielding
an omission factor for all of $sh^{-1}(I)$ of less than $r$.
By theorem 2, $sh()$ preserves density. However, $sh^{-1}()$ takes the
even numbers, which have density $\frac{1}{2}$,
to the union of \{2\} and all intervals of the form
\[ [ 2^i, \frac{3}{2} 2^i -1 ] \quad i\geq 2\]
which is a set that does not have a density.
Therefore, $sh^{-1}$ does not preserve density.

\end{document}